\documentclass[12pt]{article}
\newtheorem{thm}{Theorem}
\newtheorem{lem}{Lemma}

\newtheorem{defn}{Definition}
\newtheorem{cor}{Corollary}
\newtheorem{con}{Conjecture}
\def\ben{\begin{equation}}
\def\een{\end{equation}}

\def\bea{\begin{eqnarray}}
\def\eea{\end{eqnarray}}

\def \p{\partial}

\usepackage{chngcntr}
\counterwithout{figure}{section}
\usepackage[dvips]{graphicx}
\usepackage{pdfpages}
\usepackage{caption}
\usepackage{subcaption}

\begin{document}

\hfuzz=100pt
\title{Monotonicity of spatial critical points evolving under curvature-driven flows}
\author{G. Domokos \\
{\em Department of Mechanics, Materials, and Structures }\\
{\em Budapest University of Technology and Economics} \\
{\em M\"ugyetem rkp.3 , Budapest 1111, Hungary} \\}

\maketitle

\begin{abstract}
We describe the variation of the number $N(t)$ of spatial critical points of smooth curves
(defined as a scalar distance $r$ from a fixed origin $O$)
evolving under curvature-driven flows. In the latter, the speed $v$ in the direction
of the surface normal may only depend on the curvature $\kappa$. Under the assumption
that only generic saddle-node bifurcations occur,
we show that $N(t)$ will decrease if the 
partial derivative $v_{\kappa}$
is positive and increase  if it is negative (Theorem 1). Justification for the genericity assumption
is provided in Section 5.
For surfaces embedded in 3D, the normal speed $v$
under curvature-driven flows may only depend on the principal curvatures $\kappa, \lambda$.
Here we prove the weaker (stochastic)  Theorem 2 under the additional assumption that third-order partial derivatives can be approximated
by random variables with zero expected value and covariance.
Theorem 2 is a generalization of  a result by Kuijper and Florack for the heat equation.
We formulate a Conjecture for the case when the reference point coincides with the centre of gravity
and we motivate the Conjecture by intermediate results and an example.
Since models for collisional abrasion are governed by partial
differential equations with $v_{\kappa},v_{\lambda}>0$, our results
suggest that the decrease of the number of static equilibrium points is characteristic
of some natural processes.
\end{abstract}

%\tableofcontents

%\date{\today}       
\section{Introduction}

Curvature-driven flows are a class of nonlinear partial differential equations (PDEs),
where the time evolution of an embedded surface $\Sigma$ is defined by the speed $v$ in the direction of the surface normal
and $v$ is given as a function of the principal curvatures $\kappa, \lambda$:
\ben \label{Bloore0}
v=v(\kappa,  \lambda)
\een
In case of curves embedded in $\Re ^2,$ analogous processes are often called the \it curve shortening flows \rm \cite{Grayson} given by
\ben \label{Bloore0_2D}
v=v(\kappa)
\een
where $\kappa$ is the scalar curvature. Curvature-driven flows share
several properties with the heat equation:
\ben \label{heat}
\frac{\partial z}{\partial t}=\bigtriangleup f(x,y)
\een
(and its one-dimensional version) defining the evolution
of surface points on $f(x,y)$ not in the normal, rather in the vertical direction. Nevertheless, in the
vicinity of critical points with $\nabla f=0$, the heat equation can be regarded as an approximation to
the mean curvature flow and the latter is sometimes referred to as the geometric heat equation \cite{kimia1}. 

While locally defined, curvature-driven flows have startling global properties, e.g. they can shrink
curves and surfaces to round points (\cite{Andrews},\cite{Grayson},\cite{Huisken}). These features
made these flows powerful tools to prove topological theorems which ultimately led, via their
generalizations by Hamilton \cite{Hamilton} to Perelman's celebrated proof 
\cite{Perelman} of the Poincar\'e conjecture. The global features of curvature-driven flows are
mostly related to the monotonic change of quantities, such as the entropy associated with Gaussian curvature \cite{Chow1},
other functionals, such as the Huisken functional \cite{Huisken1} in case of the Mean Curvature flow, the number of critical points
of the curvature (used in the Curvature Scale Space model for image processing  \cite{Mokhtarian}, \cite{Mokhtarian1})
or the number $N(t)$ of spatial critical points \cite{Grayson}, closely related to the geometry of the caustic \cite{Giblin2},\cite{Giblin1}.
Beyond offering powerful tools to prove mathematical statements, curvature-driven flows also have
broad  physical applications ranging from surface
growth \cite{KPZ} through image processing \cite{Koenderink}, \cite{Mumford} to mathematical models of abrasion \cite{Bloore}, \cite{Firey}.

Our current work is primarily motivated by the latter: the monotonicity of spatial critical points
becomes an observable quantity in particle abrasion if the reference point coincides with the centre of gravity.
%>>
In this case spatial critical points appear as \it static equilibrium points, \rm i.e. points where the particle,
if it is convex, stays at rest
if supported on a horizontal surface. 
%>>
For planar curves evolving under the $v=\kappa$ flow Grayson  proved  for the case of fixed reference point (\cite{Grayson}, Lemma 1.9) the monotonic
decrease of $N(t).$ Our current work aims to prove the same for general curvature flows of type Eq.(\ref{Bloore0_2D}),
however, under the assumption that only generic, codimension-1, saddle-node bifurcations of critical points occur.
Even in case of continuous functions, monotonicity does not exclude zero gradient. However, $N(t)$ is, by definition,
an integer-valued function, so it is bound to be constant for almost all values of $t$. Monotonic decrease in this case
implies that the jumps of $N(t)$ (if and when they happen) always occur downward.

 Equation (\ref{Bloore0_2D}) describes a nonlinear partial
differential equation in rather compact notation, in section \ref{sec:2D} we will explicitly define the corresponding differential
operator (cf. equation (\ref{PDE1})). In section \ref{sec:conclusions} we will motivate our genericity assumption by pointing out that it is equivalent to admitting an arbitrarily small, additive error term in the 
aforementioned differential operator (cf. \cite{Lanczos}). On one hand, the admission of the error term  certainly simplifies the mathematical problem, thus our claim on monotonicity can be regarded
only as a first step towards the generalization of Grayson's result. On the other hand, as we will point out in section \ref{sec:conclusions}, if we regard curvature-driven flows as averaged continuum
models of abrasion processes, the genericity assumption appears plausible.
We also remark that the monotonicity of spatial critical points has also been investigated in other types of parabolic equations
in one space dimension \cite{Turyn}.

In case of surfaces evolving under curvature flows of type Eq.(\ref{Bloore0}) the monotonicity result for $N(t)$
is not true (even under our genericity assumption),
 however, Koenderink's discovery \cite{Koenderink} in 1984 about the heat equation (\ref{heat}) being the fundamental model of image representation 
sparked speculations and research about the geometrical properties of this system.
Among others, it has been early observed that creation of critical points is  \it less frequent \rm than annihilation.
As Damon \cite{Damon} notes, it was even expected that for solutions of the heat equation (\ref{heat}) no critical points would be created.
One motivation behind this belief could have been the Folklore Theorem (proven in \cite{Loog}), stating that for nonnegative,
compactly supported initial data the heat equation will eliminate all but one extremum after sufficiently long time. This theorem
certainly does not imply the monotonicity of $N(t)$ and Damon provided an example of creation of critical points which we will discuss. 
The general intuition about creation of critical points as a 'rare event' in the heat equation was later formalized by
Kuijper and Florack \cite{Kuijper} by introducing random variables as the coefficients of the Taylor series and showing that under
assuming independent, symmetric distributions the probability of creation is indeed substantially lower than that of annihilation.

We will follow the same line of thought and generalize the result of Kuijper and Florack from the heat equation (\ref{heat}) to general
curvature-driven flows given in Eq.(\ref{Bloore0}). Also, we will weaken their assumption on symmetric distribution
and only require zero expected value and covariance for the coefficients. Since principal curvatures are the only second-order surface invariants,
curvature-driven flows are the \it only \rm second-order geometric PDEs which are
invariant under coordinate transformations. This fact underscores their importance 
beyond being powerful mathematical tools, as models for physical surface evolution processes.
The general equation (\ref{Bloore0}) was first introduced by Bloore \cite{Bloore}, serving as a framework to discuss
abrasion processes.  Bloore proved some general features of Eq.(\ref{Bloore0});  he showed that the only shapes evolving under Eq.(\ref{Bloore0}) in a self-similar fashion are spheres.
He also showed that for shapes in the vicinity of the sphere, by choosing suitable scales, Eq.(\ref{Bloore0}) is equivalent to a modified version
of the heat equation (\ref{heat}) and so in this local range the qualitative features of Eq.(\ref{Bloore0}) (including the evolution of $N(t)$) are analogous to those derived for Eq.(\ref{heat}).

Our goal is to organize and extend these results to achieve a better understanding of natural abrasion processes.
First, in section \ref{sec:2D}  we will extend Grayson's Lemma for planar curves with the added genericity assumption
to more general curvature flows, including Bloore's flow.
Next, in section \ref{sec:3D} we extend our argument to the 3D case where we introduce probabilistic arguments extending the results of Kuijper  and Florack \cite{Kuijper} to general
curvature-driven flows and supporting the global,
stochastic decrease of $N(t)$ under Eq.(\ref{Bloore0}). 
Our last goal in Section \ref{sec:gravity} is to include the effect of the motion of the centre of gravity which we illustrate on an analytical example. Here again we
will argue that this effect can be modeled by a symmetric random noise
added to the system with fixed reference point. We compare the random approach with the deterministic model based on an analytic example.
In section \ref{sec:conclusions} we motivate the genericity assumption and summarize our results.

\subsection{Examples of some specific models: convergence, monotonicity and applications}
There is a rich literature on models which are special cases of Eq.(\ref{Bloore0}) and we mention some
interesting examples.
In historic order the works of Lord Rayleigh
\cite{Rayleigh1} \cite{Rayleigh2} \cite{Rayleigh3} appear to be the first where
a local, frame-invariant attrition function has been investigated; he studied
\ben \label{rayleigh}
v(\kappa,  \lambda)= c(\kappa\lambda)^{\frac{1}{4}} \,.
\een
with $c$ being constant. Rayleigh noticed that ellipsoids evolve in a self similar fashion under Eq.(\ref{rayleigh}).
Other models, also related to the Gaussian curvature $K=\kappa \lambda$
have been studied  as well, most notably Firey \cite{Firey} investigated
\ben \label{firey}
v(\kappa,  \lambda)=c\kappa\lambda
\een
(also with $c$ being constant) and proved that for any initial shape with reflection symmetry the surface converges to the sphere. 
Firey also conjectured that the symmetry requirement can be removed. This proved to be the case,  Andrews (\cite{Andrews1}, \cite{Andrews}) gave  a proof
and substantially generalized Firey's work by investigating flows under arbitrary powers of the Gaussian curvature.
The latter has broad applications ranging from abrasion processes \cite{Andrews}, \cite{Bloore}, \cite{Firey} through the affine normal flow \cite{Leicht}
to image analysis \cite{Alvarez}; a beautiful overview  of these applications is provided in \cite{Andrews2}.
Andrews \cite{Andrews} also proved some fundamental monotonicity properties of the Gauss curvature flow by defining entropy-related
functionals.
Of special interest is the mean curvature flow
\ben \label{mean}
v(\kappa,  \lambda)=b(\kappa +\lambda)
\een
(with $b$ being constant), for which Huisken \cite{Huisken} proved that it also converges to the sphere.
Huisken also proved his Monotonicity Theorem for the Huisken Functional evolving under
the Mean Curvature Flow.

Another interesting example is the PDE governing the abrasion of particles under mutual collisions.
Bloore \cite{Bloore} first discussed  general curvature-driven flows of type Eq.(\ref{Bloore0}), 
proved several fundamental properties of this equation and he also investigated the stability of the sphere by using
a power series expansion. He noted that Eq.(\ref{Bloore0}) is a rather complicated, nonlinear
geometric partial differential equation, the global solution structure of which may only
be accessible by numerical simulations.
Bloore also derived a specific model for spherical abraders with radius $R$:
\ben \label{Bloore}
v(\kappa,  \lambda)=(1+R \kappa)(1+R \lambda)
\een
and as we see, Eq.(\ref{Bloore}) is a linear combination of Eq.(\ref{firey}) and Eq.(\ref{mean}) and a constant term.
The latter may be regarded individually as a degenerate example of curvature flow:
\ben \label{Eikonal}
v \equiv 1,
\een
which is often referred to as the Eikonal equation or the parallel map and arises in the study
of wave fronts with constant speed,
satisfying Huyghen's principle. Given an initial aspherical
surface the Eikonal flow  tends to make the surface more aspherical
and to develop faces which intersect on edges \cite{DomokosSiposVarkonyi}.
In two dimensions, the ultimate shapes towards which the Eikonal flow will develop any
initial geometry are either needles with two vertices or triangles.

We also mention \it principal curvature flows \rm 
\ben \label{principal}
v(\kappa,  \lambda)=c\kappa,
\een
where the attrition speed only depends on one of the principal curvatures. This type
of evolution equation has been used in 
image processing \cite{Mumford}. The authors discuss several curvature-driven flows
serving as smoothing processes in computer vision and mention that in the ideal
case such a flow should not create new features. Among other results they
find that the number of \it umbilics \rm (i.e. points where $\kappa = \lambda$) is
certainly not monotonic under Eq.(\ref{principal}). They also point out that in two
dimensions in case of a curve evolving under a curvature-driven flow,
the number of points with zero curvature (inflection points) and
the number of extrema of the curvature is a monotonically decreasing function of time.

Curvature-driven flows have also been used in surface evolution models, most notably,
the Kardar-Parisi-Zhang
model in soft condensed matter physics describes the interface evolution via an equation for
the height function $h=h(x,y)$ as
\ben
\frac{\p h}{\p t}= \nu \nabla ^2 h  +\frac{\lambda}{2} 
(\nabla h )^2 + \eta (x,y,t) \label{KPZ}
\een 
where $\nabla$ is with respect to the flat metric on ${\bf E} ^2$ and
$\eta(x,y,t)$ is a Langevin-type  stochastic Gaussian noise term 
\cite{KPZ,Marsilli}. It was pointed out in \cite{Maritan} that this was not re-parametrization
invariant and is an approximation to a stochastic version
of the added mean curvature (\ref{mean}) and Eikonal flows (\ref{Eikonal}):
\ben \label{Mean_random}
v= -\nu H + \lambda + \eta(\varphi,\theta,t)  
\een
The first term is essentially the functional derivative 
of surface energy, i.e. a 
{\it surface tension term}  and the
second the functional derivative 
of a volume energy i.e. a {\sl pressure} term. Partially motivated by the KPZ model, in our current work we will also include stochastic terms to represent unknown quantities.
However, unlike the KPZ equation, in our model we study smooth surfaces and randomness enters the model by using random coefficients for the Fourier series expansion.

\subsection{Main results} \label{ss:mainresults}
\begin{figure}[h!]
\begin{center}
\includegraphics[scale=0.4]{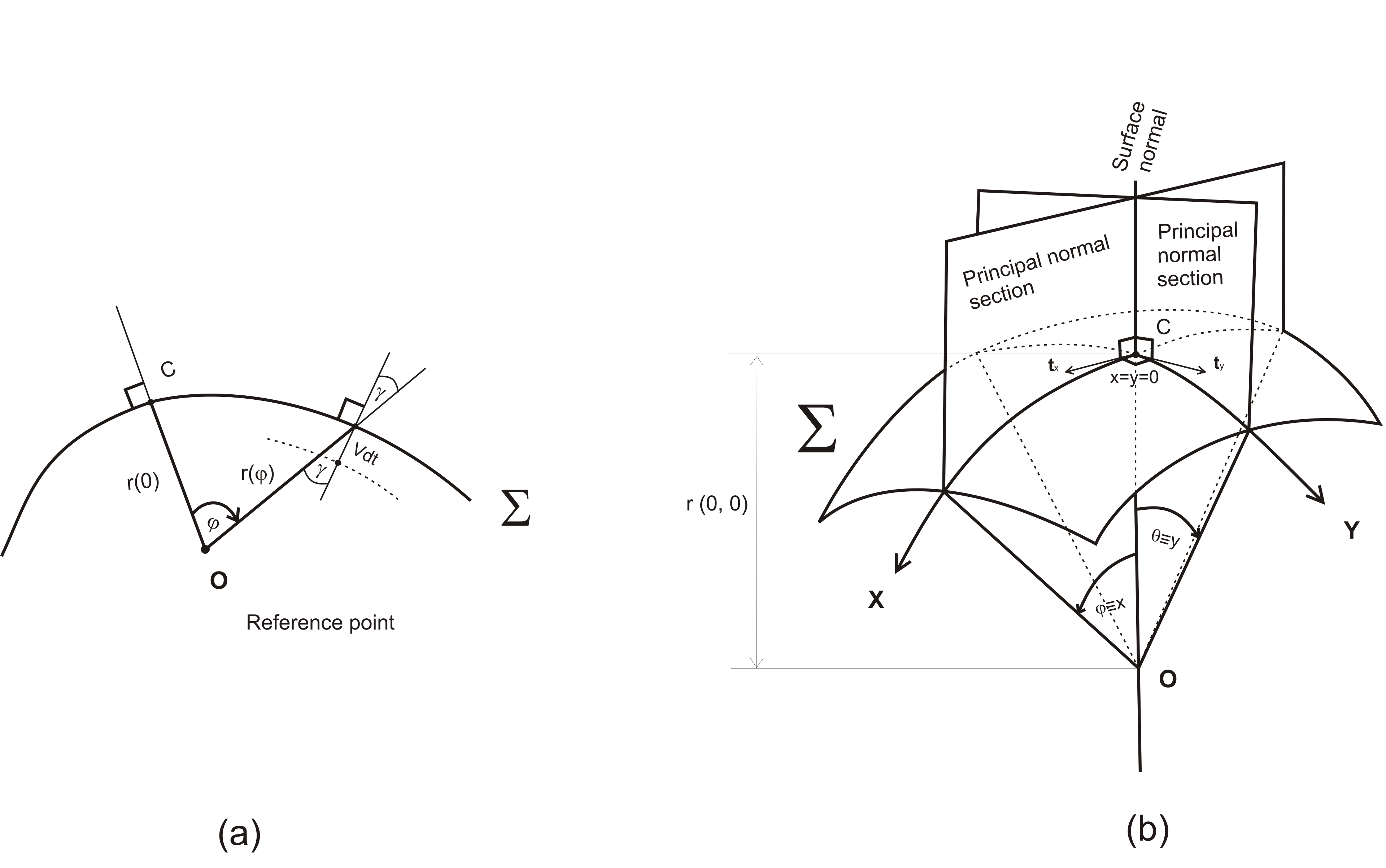}
\end{center}
\caption{\footnotesize{Local coordinates in the vicinity of the critical point $C$ of $r$: (a) 2D case, observe polar angle $\phi$
measured from $C$ and angle $\gamma$ between radial and normal direction. Cf. equation Eq.(\ref{w}), where we have $w=1/\cos \gamma$ (b) 3D case, observe polar angles $\phi \equiv x, \theta \equiv y$
with origin at $C$, measured in principal directions, defined by tangent vectors $\mathbf{t}_x,\mathbf{t}_y$. The $x,y$ coordinates at $C$ are also
tangent to curvature lines, which generally do not lie in the principal normal sections, and are not shown here.}}
\label{fig:coordinate}
\end{figure}

Our goal is to show a surprising, global feature of Eq.(\ref{Bloore0}) and Eq.(\ref{Bloore0_2D}): apparently, 
if $\Sigma$ is given as a scalar distance $r$ from a fixed reference point $O$ then the evolution $N(t)$
of the number $N$ of non-degenerate spatial critical points $C^i(t)$  $(i=1,2,\dots N(t))$ (extrema of $r(t)$ at $t=constant$) is governed primarily by the function $v$.
We will regard the  following two cases:

\begin{itemize}
	\item  the \it 2D case \rm when $\Sigma$ is a curve embedded into  $\Re ^2,$ evolving under Eq.(\ref{Bloore0_2D}),
	\item  the \it 3D case \rm
when $\Sigma$ is a surface embedded into $\Re ^3,$ evolving under Eq.(\ref{Bloore0}),

\end{itemize}
In both cases (2D and 3D) our analysis will be local, restricted to the small vicinity of the spatial critical points $C^i(t)$  $(i=1,2,\dots N(t))$ of $r(t)$. We parametrize the
spatial vicinity of $C^i$ by the polar angle $\phi$ in 2D and by two polar angles $\phi\equiv x, \theta \equiv y$ in 3D, cf. Figure \ref{fig:coordinate}. In both cases the reference point of the polar coordinate system 
is $O$ and the origin of the in-surface coordinate system on $\Sigma$ ($\phi=0; x=y=0$, respectively) coincides with the critical point $C$. Our analysis will be restricted to 
the lowest-order bifurcations of generic critical points, we make

\vspace{0.5cm}
\noindent \bf Assumption 1:\rm 
\it We assume that as the scalar distance function $r(t)$ evolves under Eq.(\ref{Bloore0}) or Eq.(\ref{Bloore0_2D}), spatial
critical points $C^i(t)$ of $r(t)$ will undergo only generic, codimension-1, saddle-node bifurcations. \rm

\vspace{0.5cm}

\noindent  Assumption 1 is of central importance in the paper and all results in sections \ref{sec:2D},
\ref{sec:3D} and \ref{sec:gravity} are based on this hypothesis. In section \ref{sec:conclusions} we motivate Assumption 1 by pointing out
that it is equivalent to the admission of  arbitrarily small perturbation of the differential operators corresponding to Eq.(\ref{Bloore0}) and
Eq.(\ref{Bloore0_2D}). Our main motivation is to use these equations as models for physical processes because there is strong evidence \cite{jerolmack}
based on laboratory experiments that they indeed serve as good models. In this context we
expect that the information lost due to arbitrarily small perturbations of the operators is not relevant
for these processes, and we regard Assumption 1 as sufficient for our purpose. We also note that systems with codimension-1, saddle-node bifurcations form an open-dense
set in the space of smooth 1-parameter dynamical systems.

In section \ref{sec:2D}, by relying on Assumption 1 and drawing on some simple facts from bifurcation theory, we will prove 
\begin{thm}\label{thm1}
In the 2D case, under Assumption 1, $N(t)$ is monotonically decreasing under Eq.(\ref{Bloore0_2D}) if and only if $v_{\kappa} > 0$,
\end{thm}
\noindent where subscript denotes partial differentiation.

In the 3D case, the variation of $N(t)$ is controlled
by two additive terms, one of which is analogous to the planar case and
has the opposite sign of $v_{\kappa}$, the sign and magnitude
of the second term depends on the third-order partial derivatives ${r_{xxy}},{r_{yyy}}$ which are, in general, unknown. The same difficulty arises in the context
of the heat equation. Kuijper and Florack \cite{Kuijper} introduced a probabilistic assumption 
and we follow this by introducing the slightly weaker

\vspace{0.5cm}
\noindent
\bf Assumption 2 \rm \it At the critical point  $C$, we replace ${r_{xxy}},{r_{yyy}}$ by the 
random variables ${\bar r_{xxy}},{\bar r_{yyy}}$ with $E(\bar r_{xxy})=E(\bar r_{yyy})=\sigma({\bar r_{xxy}},{\bar r_{yyy}})=0.$ \rm

\vspace{0.5cm}

\noindent  Under Assumption 2, $N(t)$ is replaced by a stochastic process with expected value $\bar N(t)$ and 
in section \ref{sec:3D} we will prove
\begin{thm}\label{thm2}
In the 3D case, under Assumptions 1 and 2, $\bar N(t)$ is monotonically decreasing under Eq.(\ref{Bloore0}) if and only if $v_{\kappa},v_{\lambda} > 0.$ 
\end{thm}

We remark that Kuijper and Florack also (tacitly) assumed symmetric distributions for ${\bar r_{xxy}},{\bar r_{yyy}}$, however, this is not
needed for our current purpose and, unlike zero expected value, it would be hard to justify. We also remark
that, similarly to the results of Kuijper and Florack, neither the exact form of the distributions for $\bar r_{xxy}, \bar r_{yyy}$ nor their
respective  individual variances enter
into our formulae and thus it is sufficient to make an assumption about their expected values and covariance being zero.
In section \ref{sec:3D} we will show that an independent and apparently weaker probabilistic assumption
(which we call Assumption 2-A) is a suitable alternative to Assumption 2 and leads also to Theorem 2.

So far we discussed problems related to a fixed reference point $O$. In physical applications, such as abrasion processes, designating
a fixed reference point is often very problematic, so we also regard the case of moving reference.
If $O$ is identical to the centre of gravity $G$ then
stationary points $C^i$ of $r$ correspond to \it static equilibria. \rm  In this case, based alone on the local PDE (\ref{Bloore0})
we do not have immediate information on the evolution of $N(t)$ since the location of $G$ is, in general, not fixed with
respect to the body. However, in Section \ref{sec:gravity} we will formulate Conjecture \ref{con:gravity} stating that $\bar N(t)$ remains monotonically
decreasing even in this case. To motivate the Conjecture, we demonstrate special cases (subsection \ref{ss:symmetry}) and a stochastic assumption
(subsection \ref{ss:stochastic}) under which the conjecture holds.

The two Theorems and the Conjecture are in increasing order of physical applicability, however (due to heuristic assumptions
about randomness of unknown quantities) in decreasing order of mathematical strength.
The properties of Eq.(\ref{Bloore0}) and Eq.(\ref{Bloore0_2D}) stated by the Theorems are of special interest because,
as we can see, in many of the listed examples we
have $v_{\kappa},v_{\lambda}>0,$ so our arguments predict a stochastically decreasing
trend for the number of static equilibria for convex bodies evolving under any of these models.

\section{Evolution of critical points on planar curves} \label{sec:2D}
\subsection{Evolution assuming generic bifurcations}
Our goal is to prove Theorem 1.
We regard a closed curve $\Sigma (t) $ embedded in $\Re^2$, evolving  under Eq.(\ref{Bloore0_2D}).
The curve is given by the  $C^\infty$-smooth distance function $r(t)$, measured from a
 fixed origin $O$ in the interior of $\Sigma(t)$.
All our analysis will be restricted to the small vicinity of spatial critical points $C^i(t) \in \Sigma$, $(i=1,2,N(t))$ of $r(t)$, and at such points
 $r$ can be always locally parametrized by the polar angle $\phi$ and the curvature
$\kappa$ can be expressed \cite{Oxford} as
\ben \label{curvature}
\kappa(r,r_{\phi},r_{\phi \phi}) = \frac{r^2+2r_{\phi}^2-rr_{\phi\phi}}{(r^2+r_{\phi}^2)^{\frac{3}{2}}}.
\een
(Note that the above polar curvature formula was already known to Newton \cite{Whiteside}.)
Now we may write the speed in the normal direction, instead of Eq.(\ref{Bloore0_2D}) as:
\ben \label{Planarbloore1}
v=v(r,r_{\phi},r_{\phi \phi}).
\een
Note that Eq.(\ref{Bloore0_2D}) is a special case of Eq.(\ref{Planarbloore1}) and based on our proof
the result can be also interpreted in this more general setting.
The radial speed is given as
\ben \label{PDE1}
r_t=-v(r,r_{\phi},r_{\phi \phi})w(r,r_{\phi})
\een
where subscript $t$ denotes partial derivative with respect to time $t$ and 
\ben \label{w}
w=\sqrt{\frac{r^2+r_{\phi}^2}{r^2}}
\een
is accounting for the tangential movement of the points 
(cf. Figure \ref{fig:coordinate}(a), where we can see that $w=1/\cos \gamma$).
Under Assumption 1, critical points appear or disappear only at codimension-1, saddle-node
bifurcations where $r(\phi,t)$ is locally diffeomorphic to 
\ben \label{saddlenode0}
r=\phi^3 \pm t \phi,
\een 
and for the derivatives we have
\ben \label{saddlenode}
r_{\phi}=r_{\phi\phi}=0, \qquad r_{\phi\phi\phi}\not=0.
\een
After differentiating Eq.(\ref{PDE1}) with respect to  $\phi$
and substituting Eq.(\ref{saddlenode}) we get
\ben \label{c0}
r_{t\phi}=-v_{\phi}=-v_{r_{\phi\phi}}r_{\phi\phi\phi}.
\een
At such a saddle-node a pair of critical points is either created or annihilated.
Which of these two possibilities is realized depends
on the relative sign of the cubic term $r_{\phi\phi\phi}$ and the emerging linear term $r_\phi$.
Since the latter is zero at the bifurcation, its time derivative $r_{t\phi}$ will determine
its sign after passing the bifurcation. Alternatively, we see that.
the $\pm$ sign in Eq.(\ref{saddlenode0}) determines whether critical points are created
or annihilated. We introduce the \em annihilation
indicator \rm
\ben \label{omega0}
\Omega \in \{-1,0,1\}
\een
 such that $\Omega=1$ corresponds to annihilation of critical points and $\Omega =-1$
to the creation of new ones. 
 In case of  generic saddle nodes on smooth, planar curves we have
\ben \label{annih}
\Omega=\mbox{sgn}\left(\frac{r_{t\phi}}{r_{\phi\phi\phi}}\right),
\een
and, based on Eq.(\ref{c0}) we have
\ben \label{annih1}
\Omega =\mbox{sgn}\left( -v_{r_{\phi\phi}}\right).
\een
We can now re-write Eq.(\ref{annih1}) in terms of the invariant curvature $\kappa$ given by Eq.(\ref{curvature}).
Based on the latter equation, at the  saddle node we have
\ben \label{kappa1}
\kappa = \frac{1}{r}, \qquad \kappa_{r_{\phi\phi}}=-\frac{1}{r^2}=-\kappa ^2
\een
and condition Eq.(\ref{annih1}) can be re-written as
\ben \label{kappa2}
\Omega=\mbox{sgn}\left(-v_{r_{\phi\phi}}\right)=\mbox{sgn}\left(-v_{\kappa}\kappa_{r_{\phi\phi}}\right)=\mbox{sgn}\left(\kappa ^2 v_{\kappa}\right)=\mbox{sgn}\left(v_{\kappa}\right).
\een
Using this concept we can write the following condition:
a new pair of critical points will be created or annihilated at a saddle node if
\ben \label{c1}
v_{\kappa} < 0, \qquad v_{\kappa} > 0,
\een
respectively. Thus we have proven Theorem 1.

\noindent We also mention that in the planar ($\lambda=0$) version of the specific Bloore equation (\ref{Bloore}) we have
\ben \label{Omega_bloore}
v_{\kappa}=R, \qquad \Omega_{Bloore} = \mbox{sgn}\left(\kappa^2 R\right)=1,
\een
so $N(t)$ will be monotonically decreasing under the planar version of Eq.(\ref{Bloore}).

\subsection{The pitchfork: an illustration} \label{ss:pitchfork}

So far we looked at generic, codimension-1 saddle nodes. Here we
would like to indicate that our argument can be carried over to pitchforks
which are degenerate, however, still structurally stable.
Such bifurcations could arise, for example, as a result of reflection symmetry.
In the planar case, if the symmetry axis passes through 
the investigated critical point of $r(\phi)$ then this critical point is
degenerate ($r_{\phi}=r_{\phi\phi}=0$) and new critical points
are created/absorbed not in 
saddle-nodes (described in Eq.(\ref{saddlenode})), rather at pitchforks 
where the first $n=3$ terms in the Taylor series expansion $T(r(\phi))$
of $r(\phi)$ vanish and the $(n+1)=4th$ term does not vanish:
\ben \label{pitchfork}
r_{\phi}=r_{\phi\phi}=r_{\phi\phi\phi}=0, \qquad r_{\phi\phi\phi\phi}\not=0.
\een
In this case, based on Eq.(\ref{c0}) we have $r_{t\phi}=0$
and the annihilation indicator is defined as
\ben \label{annih_pitchfork}
\Omega=\mbox{sgn}\left(\frac{r_{t\phi\phi}}{r_{\phi \phi\phi\phi}}\right).
\een
After differentiating Eq.(\ref{c0}) 
 and using Eq.(\ref{pitchfork}) we get
\ben \label{c00}
r_{t\phi\phi}=-v_{\phi\phi}=-v_{r_{\phi\phi}}r_{\phi\phi\phi\phi},
\een
and we arrive at
equation (\ref{kappa2}). So Theorem 1 remains valid in this case.
\begin{figure}[h!]
\begin{center}
\includegraphics[scale=0.5]{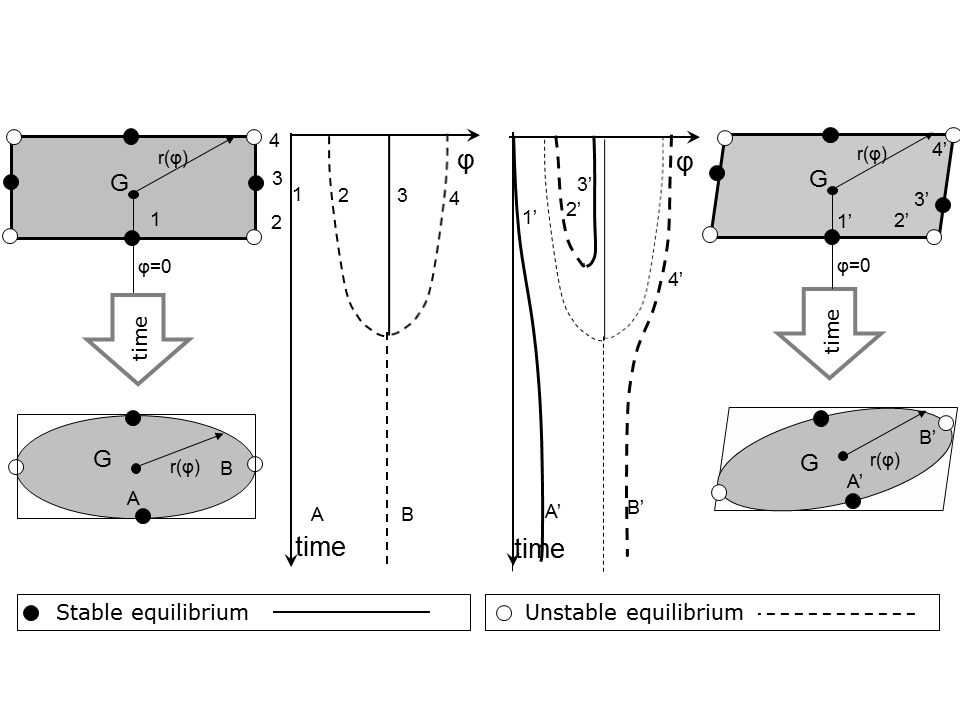}
\end{center}
\caption{\footnotesize{Qualitative evolution of stationary points under curvature-driven abrasion.
Rectangular block (upper left corner) with 8 critical points is abraded into
ellipsoidal shape with 4 critical points. Bifurcation diagram of critical points
on the right, dashed line correspond to maxima, solid lines to minima.
Small distortion of rectangle into a parallelogram (upper right corner) breaks
the reflection symmetries and results in unfolding
the pitchfork into a saddle-node and a nearby continuous path.}}
\label{fig:1}
\end{figure}

Such a pitchfork appears if a rectangular block is being abraded by curvature-driven
flow $v=c\kappa$ to a roughly ellipsoidal shape.  Even without detailed computations this
process is intuitively clear. Initially the  distance function $r(\phi)$ for  the rectangle has $N=8$ critical points
(4 maxima and 4 minima), the final, smooth shape will have $N=4$ critical points
(2 maxima and 2 minima). We illustrate the corresponding bifurcation diagram
in Figure \ref{fig:1}, left side. If we tilt the vertical sides of the rectangular block by the same
small angle then we obtain an imperfect problem represented by
a parallelogram which breaks both reflection symmetries
of the original problem (cf. Figure \ref{fig:1}, right side). The parallelogram will also be abraded to an ellipsoid-like
shape and the numbers for the critical points agree in the symmetric and imperfect
problem both for the initial and the final shapes. However, the evolution
is different: in the symmetrical case we observe pitchforks while in the
imperfect case we see saddle-nodes. In both cases the number $N$ of critical
points decreases.

\section{Evolution of critical points on surfaces}\label{sec:3D}
Our goal is to prove Theorem 2.
We now regard a closed, $C^\infty$-smooth surface $\Sigma$ embedded in $\Re^3$ and a fixed reference
point $O$ in the interior of $\Sigma$. The scalar distance between $O$ and $\Sigma$ will be denoted by $r$
which is a function in two variables.
At each point $P$ on $\Sigma$ (except for umbilic points) we can define a local, orthogonal, in-surface $xy$ coordinate system by the normal
planar sections of $\Sigma$ in the directions of the principal curvatures, see Figure \ref{fig:coordinate}. Our local Cartesian coordinates are tangent to the
in-surface \em curvature lines. \rm  At umbilic points the directions of
principal curvatures are not defined since all curvatures are equal.
Stationary points $C^i, (i=1,2,\dots N(t))$ 
of $r$ are characterized by 
\ben \label{equilibrium_3D}
r_x=r_y=0, \quad r_{xy}=0, \quad r_{xx}\neq 0, r_{yy}\neq 0.
\een
(The equality  $r_{xy}=0$ is due to our special choice of coordinates since the Hessian is diagonal in the
frame aligned with the principal curvatures.)
In the vicinity of the critical point we parametrize our coordinate lines $x,y$ by the  spherical polar angles $\varphi,\theta$
with origin at $O$. Since the latter is on the surface normal at $C$,
by suitable shift of the coordinates  we can always obtain 
\ben \label{xy}
x \equiv \varphi, \qquad y \equiv \theta,
\een 
and at $C$ we have
\ben \label{origin}
x=y=0.
\een
Now the values of both principal curvatures
 $\kappa, \lambda$ may be computed  at $(x,y)=(0,0)$ by substituting Eq.(\ref{xy}) into formula (\ref{curvature}) 
containing the expression of curvature in planar polar coordinates:
\bea \label{curvature_3D}
\lambda^C(r,r_{x},r_{xx}) & = & \frac{r^2+2r_{x}^2-rr_{xx}}{(r^2+r_{x}^2)^{\frac{3}{2}}} \\
\label{curvature_3D_1}
\kappa^C(r,r_{y},r_{yy}) & =  & \frac{r^2+2r_{y}^2-rr_{yy}}{(r^2+r_{y}^2)^{\frac{3}{2}}}.
\eea
We note that the principal curvatures (\ref{curvature_3D}-\ref{curvature_3D_1}) are not  invariant under
arbitrary changes in the $x,y$ coordinates as these are polar angles and they rescale simultaneously
with the radius $r$ if the point of reference $O$ is being moved on the surface normal.
We also note that formulae (\ref{curvature_3D}-\ref{curvature_3D_1}) provide the values of the principal curvatures only at $C$ (indicated by the superscript $C$),
however, the partial derivatives of $\lambda, \kappa$ with respect to $x$ and $y$ may not, in general, be obtained by formally
differentiating Eqs. (\ref{curvature_3D}-\ref{curvature_3D_1}) since curvature lines are, in general, non-planar, i.e. they do not coincide with
our coordinate lines.

Now we extend $\Sigma$ to be a member of
a generic, one-parameter family ${\Sigma}(t)$ of smooth surfaces. In such a family, at isolated values of $t,$ codimension-1, generic
saddle node bifurcation points may occur
where the Hessian of $r(x,y)$ becomes singular and under Assumption 1, we can restrict our analysis
to the vicinity of such points. Since we use the $[x,y]$ directions of principal curvatures as local coordinates, the
Hessian is diagonal and thus its 
singularity requires that one of the second-order derivatives $r_{xx},r_{yy}$ vanishes \cite{Poston}
so, based on Eq.(\ref{equilibrium_3D}), here we have
\ben \label{saddlenode2}
r_x=r_y=r_{xy}=r_{yy}=0, \quad r_{xx}\not=0, r_{yyy}\not=0, r_{xxy}\not=0, r_{yyx}\not=0, r_{xxx}\not=0
\een
At such a codimension-1, generic saddle-node the principal curvatures can be written based on 
Eqs.(\ref{curvature_3D}-\ref{curvature_3D_1}) and Eq. (\ref{saddlenode2}) as
\bea
\lambda^C & = & \frac{r-r_{xx}}{r^2} \\
\label{curvature_3D_1a}
\kappa^C & =  & \frac{1}{r},
\eea
so we always have $\lambda^C \not=\kappa^C$ and thus this cannot be an umbilic point of $\Sigma$, i.e.
 our $[x,y]$ system is always well-defined at codimension-1, generic saddle nodes.
At such points, a pair of critical points is either created
or annihilated and based on  Eq.(\ref{saddlenode2}), $r(x,y)$ can be approximated up to third-order by
\ben \label{3rdorder1}
r(x,y)\approx r+\frac{r_{xx}}{2}x^2+\frac{r_{yyy}}{6}y^3+\frac{r_{xxy}}{2}x^2y+\frac{r_{yyx}}{2}y^2x+\frac{r_{xxx}}{6}x^3,
\een
where the functions $r,r_{xx},r_{yyy},r_{xxy},r_{yyx},r_{xxx}$
on the right hand side are evaluated at $(x,y)=(0,0)$.
We note that higher order terms do not influence either the curvatures or their first
derivatives, so this approximation is sufficient for our current purpose.

Next we discuss how the generic, codimension-1 saddle-node evolves under curvature-driven flows. We let both
the reference point $C$ and the coordinate lines $[x,y]$ be propagated by the flow which is at the saddle-node locally equivalent to
a normal flow. The evolution equation
defining ${\Sigma}(t)$ at the saddle-node is given by the (scalar) speed $v$ in the direction of the inward surface
normal. Since the principal curvatures $\kappa, \lambda$ do not depend on partial derivatives beyond second order, based on Eq.(\ref{Bloore0})we can write
\ben
v=v(r,r_x,r_y,r_{xx}, r_{xy}, r_{yy}),
\een
and we note that, based on Eqs. (\ref{curvature_3D}-\ref{curvature_3D_1}), at the saddle-node the principal curvatures do not depend on the mixed derivative $r_{xy}$.
The evolution equation for $\Sigma$ in the radial direction can be written, analogously to Eq.(\ref{PDE1}) as
\ben \label{PDE_3D}
r_t=-v(\kappa, \lambda)w(f(y),r,r_y,r_x)=-v(r, r_x,r_y, r_{xx}, r_{xy}, r_{yy})w(f(y),r,r_y,r_x)
\een 
where $w= \sqrt{(r^2_x + r^2_y \cos^2 y + r^2 \cos^2 y)/(r^2 \cos^2y)}$.
If we differentiate Eq.(\ref{PDE_3D}) at $(x,y)=(0,0)$ we get:
\ben \label{3Dc0}
r_{ty}=-v_{\lambda}\lambda _y-v_{\kappa}\kappa _y
\een
As in the planar case,
a new pair of critical points will appear or disappear at a saddle node if
\ben \label{3Dc1}
\frac{r_{ty}}{r_{yyy}} < 0, \qquad \frac{r_{ty}}{r_{yyy}} > 0.
\een
The annihilation indicator $\Omega$ defined in Eq.(\ref{annih}) can be expressed from Eq.(\ref{3Dc0}) as:
\ben \label{3Dc3}
\Omega=\mbox{sgn}\left(\frac{-v_{\lambda}\lambda _y-v_{\kappa}\kappa _y}{r_{yyy}}\right).
\een
To obtain the partial derivatives $\kappa _y, \lambda _y$ 
at $(x,y)=(0,0)$
we used Eq.(\ref{3rdorder1}) and with the aid of  direct computation with Maple 16 we arrived at
\bea \label{ly}
\kappa_y & = & -\frac{r_{yyy}}{r^2} \\
\label{ky}
\lambda_y & = & -\frac{r_{xxy}}{r^2} .
\eea
We remark that we arrive at the same formulae if we differentiate Eqs. (\ref{curvature_3D}-\ref{curvature_3D_1}),
however, the explanation appears to be not obvious. We also remark that the formulae for $\lambda _x, \kappa _x$
may not be obtained by formally differentiating Eqs. (\ref{curvature_3D}-\ref{curvature_3D_1}), however, 
we do not need these quantities for our computations. The difference between the behavior of the $x$ and $y$ derivatives
is rooted in the asymmetry of the truncated Taylor series (\ref{3rdorder1}).

Substituting Eqs. (\ref{ly}-\ref{ky}) into Eq.(\ref{3Dc3}) yields
\ben \label{alpha2}
\Omega= \mbox{sgn}\left(v_{\kappa} + \frac{r_{xxy}}{r_{yyy}}v_{\lambda}\right).
\een
As we can observe, in the 3D case the critical point absorbing/generating property of the
PDE can, in general,  not be decided based alone on the formula for the attrition speed $v$.
Apparently, this problem has been deeply investigated in the context of the
closely related heat equation (\ref{heat}) which can be also written as
\ben \label{heat1}
r_t=-v=r_{xx}+r_{yy}.
\een
Damon \cite{Damon} provided an explicit example for the creation of critical points,
in our notation his example is:
\ben
r(x,y,t)=y^3-6ty-6yx^2+x^2+2t
\een
which satisfies the heat equation. As we can see, at the critical point $r(0,0,0)$ we have
\ben
r_x = r_y=r_{xy}=r_{yy}=0, \quad r_{yyy}=6, \quad r_{xxy}=-12,
\een
and so, based on Eq.(\ref{3Dc1}) we have a creation.

While apparently the 3D case admits both annihilations and creations, the latter have been observed to occur less frequently
(\cite{Damon} \cite{Loog}). This observation certainly does not imply the monotonicity of $N(t)$, nevertheless, one would expect that in some averaged sense $N(t)$ will decrease.
We will formalize this intuition by introducing random variables.
First we observe that we can multiply the right hand side of  Eq.(\ref{alpha2}) by $r_{yyy}^2$ without changing the sign of the left hand side, so we get
\ben \label{alpha2a}
\Omega= \mbox{sgn}\left(v_{\kappa}r_{yyy}^2 + {r_{xxy}}{r_{yyy}}v_{\lambda}\right).
\een

Following the ideas of Kuijper and Florack \cite{Kuijper}
we regard ${r_{xxy}},{r_{yyy}}$ as random variables with zero expected values and covariance, however, we do not require
that the distributions should by symmetrical. We remark that assuming zero expected value can be geometrically motivated
and this is most easily seen if $\Sigma$ is convex since in this case we can globally extend our $[xy]$ spherical polar coordinates.
If we choose $y$ to be the lateral angle ($y \in [0,2\pi)$) of the spherical polar system then we have periodicity in this
variable. In particular, the average values of partial derivatives in the $y$ direction will be zero.
Since both $r_{xxy}$ and $r_{yyy}$ fall into this category,
it is plausible to adopt Assumption 2 (subsection \ref{ss:mainresults}) replacing them with random variables with zero expected value and covariance.
Now $\Omega$ is also replaced by the random variable $\bar \Omega$ and based on Eq.(\ref{alpha2a}) we can write
\ben \label{alpha2b}
E(\bar \Omega)=E\left( \mbox{sgn}\left(v_{\kappa}\bar r_{yyy}^2 + {\bar r_{xxy}}{\bar r_{yyy}}v_{\lambda}\right)\right).
\een
Since $E(\bar r_{yyy}^2) >0$, under Assumption 2 we have
\ben
\mbox{sgn}(E(\bar \Omega)) =\mbox{sgn}(v_{\kappa}).
\een
and so we have proven Theorem 2 under Assumption 2.
If, in addition to the latter we also assume that ${\bar r_{xxy}},{\bar r_{yyy}}$ are symmetric with respect to zero
and in Eq.(\ref{Bloore0}) we assume $v_{\kappa}=v_{\lambda}$ then we have
\ben
E(\bar \Omega)=0.5, \qquad P(\bar \Omega =1)=75\%
\een 
and this was pointed out by Kuijper and Florack (\cite{Kuijper}) in case of the heat equation (\ref{heat1}). We note that the assumption on
the symmetry of the random variables is not easy to justify and so these numerical values may just serve as an illustration.

We remark that by introducing random variables  in Eq.(\ref{alpha2b}), the function $N(t)$ is replaced by a one-dimensional, continuous time
 Markov-chain the state space of which are, based on Assumption 1, the even natural numbers
and subsequent states always differ by $\Delta N=\pm 2$ (in case of additional symmetry, subsequent states may differ by other even numbers).
We can define this process in two stages: first we define a discrete time Markov chain with
transition probabilities
\ben \label{transition}
P(\Delta N= 2i) = P(\bar \Omega =-i), \quad i = \pm 1.
\een
In the second stage the random variables associated with the holding times in each state may be defined independently which we do not specify here.
Theorem 2 can be interpreted in this setting by the claim that the discrete-time Markov process (defined in the first stage in Eq.(\ref{transition}))
is asymmetric, the direction of the drift is opposite to the sign of $v_{\kappa}$.

We will next show that an alternative, independent random assumption also leads to the same conclusion.
First we observe that under Assumption 1, at a generic codimension-1, saddle node bifurcation $\bf C^{\star} \rm $ of $r$ we have a generic (Morse) critical point
$C_y^{\star}$ for $r_y(x,y)$(cf. also equation (\ref{saddlenode2})). We will prove

\begin{lem} \label{lemma:elliptic}
If the Morse critical point $C_y^{\star}$ of $r_y$ coinciding with the saddle-node $\bf C^{\star} \rm $ of $r$ 
is elliptic then $N(t)$ is decreasing at $\bf C^{\star}\rm $.
\end{lem}
\noindent For the determinant Det${\bf H}$ of the Hessian of $r_y(x,y)$ we have
\ben
\mbox{Det}{\bf {H}} \leq r_{xxy}r_{yyy},
\een
so based on Eq.(\ref{alpha2a}) we can write
\ben \label{alpha2c}
\Omega \geq \mbox{sgn}\left(v_{\kappa}r_{yyy}^2 + \mbox{Det}{\bf H}v_{\lambda}\right).
\een
If the critical point of $r_y(x,y)$ is elliptic then Det${\bf H}>0$, thus
we have proven Lemma \ref{lemma:elliptic}. Next we introduce an alternative to Assumption 2 in subsection \ref{ss:mainresults}:

\noindent \bf Assumption 2-A \rm :
\it We assume that $\Sigma$ is a topological sphere and the critical point $C_y^{\star}$ of $r_y$, coinciding with the saddle-node $\bf C ^{\star}$ of $r,$ is selected by uniform random choice from among the
Morse critical points $C_y^i$ of $r_y$.\rm

\noindent Now we proceed to prove Theorem \ref{thm2} under Assumption 2-A.
The Poincar\'e-Hopf formula  relates the numbers of minima, maxima and saddles ($S,U,H$, respectively)
on a topological sphere as
\begin{equation} \label{Poincare}
\chi=S+U-H=2
\end{equation}
and $\chi$ is also called the Euler characteristic of the surface.
We observe that $r_y(x,y)$ is the partial derivative of a periodic function, so its average
value is zero, therefore its \emph{global} extrema cannot be zero. Therefore the
numbers $S,U$ of maxima and minima from among which the stationary point is picked have to
be reduced by 1 to obtain the relevant numbers 
\ben  \label{relevant}
\bar S=S-1, \quad \bar U=U-1, \quad  \bar H = H
\een
and the total number of relevant critical points is $n=\bar S + \bar U + \bar H= S+U+H-2$.
We denote the relevant critical points of $r_y(x,y)$ by
$C^i_y, i=1,2,\dots n$ and the Hessian determinants associated with these points by Det$_i$.
Using Assumption 2-A, Eq.(\ref{alpha2c}) is replaced by
\ben \label{alpha3}
\bar \Omega= \mbox{sgn}\left(v_{\kappa}r_{yyy}^2 + \xi v_{\lambda}\right),
\een
where $\xi$ is a random variable with $P(\xi=\mbox{Det}_i)=1/n$ and $\bar \Omega$ is a random variable and we seek the sign of its expected value.
Based on Eq.(\ref{Poincare}) and Eq.(\ref{relevant}) we can see that 
\ben
\bar S + \bar U = \bar H
\een
i.e. the number of relevant elliptic and relevant hyperbolic points is equal. Using Lemma (\ref{lemma:elliptic}),
the random variable $\chi$ will thus assume positive and negative values with equal relative frequency
so we have
\ben
E(\mbox{sgn}(\xi))=0,
\een
and thus we have
\ben \label{result2}
\mbox{sgn}(E(\bar \Omega))= \mbox{sgn}(v_{\kappa}).
\een
So we have proven Theorem 2  also under Assumption 1 and Assumption 2-A.

We remark that for closed, connected  2D manifolds the Euler characteristic $\chi$ may
be much lower than 2 (e.g. for the simple torus we have $\chi=0$).
We also remark that in the specific Bloore flows Eq.(\ref{Bloore}), based on Eq.(\ref{curvature}) and Eq.(\ref{result2}) we have
\ben
\mbox{sgn}(E(\bar \Omega))=\mbox{sgn}((\kappa ^2 R(R\lambda+1))=1,
\een
so we see that the expected value $\bar N(t)$ of the stochastic process is decreasing under Eq.(\ref{Bloore}).

\section{Centre of gravity: moving reference point}\label{sec:gravity}
So far we dealt with critical points with respect to a fixed reference point $O$; for brevity  we will
call this the \it fixed model. \rm If we regard curvature-driven flows
as models for abrasion processes then the notion of fixed reference has to be carefully examined.
The process itself is invariant under any change of coordinates and as the surface evolves, it is
not trivial to define any fixed point in a physically meaningful manner.
One possible choice for a physically objective reference is the \it ultimate point \rm $U$ which is abraded last at $t=t_u$.
One drawback of this choice is that the distance from this point does not appear to have any special physical meaning.

From the point of view of applications it
is of interest to study the case when $O\equiv G$ i.e. when the reference point is the centre of gravity.
In this case, critical points of the distance function
coincide with \it static equilibrium points \rm of the object, i.e. surface points on which the convex objects under vertical gravity
 can be balanced when
supported on a horizontal surface. (In case of nonconvex minima the support surface has to have sufficient curvature.)
 We will refer to this approach as the \it centroid model \rm and we formulate

\begin{con} \label{con:gravity}
In the centroid model $\bar N(t)$ is decreasing monotonically.
\end{con}

In the rest of section \ref{sec:gravity} we try to motivate Conjecture \ref{con:gravity} by pointing out special cases (subsection \ref{ss:symmetry}) and an additional stochastic Assumption (subsection \ref{ss:stochastic})
under which the Conjecture can be proven. We illustrate the latter with an analytical example (subsection \ref{ss:example}).

\subsection{Symmetry}\label{ss:symmetry}

Let the joint symmetry group (isometry group) of the surface $r(\phi, \theta)$ and the distribution of material density be denoted
by $\Gamma$. By definition, the location of $G$ is invariant under $\Gamma$. If the invariant subspace of $\Gamma$
is a single point  then we say that  the object has a \it unique centre of symmetry. \rm In this case $G$ will be fixed under any local PDE, so our results for the fixed model apply directly:

\begin{cor} \label{cor:symmetry}
to Theorem \ref{thm2}. In the 3D case, under Assumptions 1 and 2, if $r$ is measured from the centre of gravity $G$ and the object has a unique centre of
symmetry , then 
$\bar N(t)$ is monotonically decreasing under Eq.(\ref{Bloore0}) if and only if $v_{\kappa},v_{\lambda} > 0$ 
\end{cor}

We also note that  $C^{\infty}$-small, generic perturbations of symmetric objects will preserve the number of critical points, so as long
as the object remains in the $C^{\infty}$-vicinity of an object with unique centre of symmetry, $\bar N(t)$ will remain monotonic.
If, however, this is not the case then the motion of $G$ is not controlled by the local properties 
of the surface in the vicinity of its critical points. As the body's shape changes by losing material,
$G$ moves with respect to fixed points in the body and integration over the evolving shape
is needed to determine its position. In this case the centroid model and the fixed model will, in general, predict different evolution for $N(t)$.

\subsection{A stochastic assumption for the centroid model}\label{ss:stochastic}

Our goal is to approximate the centroid model by a fixed model with added, symmetric noise, the \it stochastic centroid model \rm
which is defined by

\noindent
\bf Assumption 3 \rm \it We regard the motion of $G$ (determined by integration over the evolving shape) to be independent of the local properties of the surface in the vicinity
of the equilibrium point. We model the effect of the motion of $G$ on the value of $\Omega$ by an added, symmetric random noise $\eta$ uniform on $[-d,d].$\rm
\noindent

First we show that under Assumption 3, Conjecture \ref{con:gravity} holds.
Under Assumption 3, instead of Eq.(\ref{alpha3}) we write
\ben \label{alpha4}
\bar \Omega= \mbox{sgn}\left(v_{\kappa}r^2_{yyy} + \xi v_{\lambda}+\eta \right) ,
\een
In the planar case we have
\ben \label{alpha5}
\bar \Omega= \mbox{sgn}\left(v_{\kappa}  +\eta \right).
\een
Since $E(\eta)=0$, analogously to Eq.(\ref{result2}), in Eq.(\ref{alpha4}),Eq.(\ref{alpha5})  we get
\ben \label{result3}
\mbox{sgn}(E(\bar \Omega))= \mbox{sgn}(v_{\kappa}).
\een
So we see that Theorem \ref{thm2} remains valid if we add Assumption 3 to Assumptions 1 and 2.
To motivate Assumption 3, we will illustrate the difference between the fixed and centroid models
and their approximation by the stochastic centroid model on a planar, analytical example under the Eikonal flow.
(In case of a 3D example under genuine curvature flows, undoubtedly closer to the main objective of the paper,
only numerical treatment appeared feasible. We decided against the latter because it is well known that  even arbitrarily fine discretizations can yield spurious solutions
\cite{Peitgen}, recently it was even demonstrated \cite{Langi} that they also yield additional, spurious critical points
on curves and surfaces.)

\subsection{The Eikonal flow}

The degenerate case of the parallel (Eikonal) flow $v \equiv 1$ (\ref{Eikonal}) is, strictly speaking,
 certainly  \it not \rm a curvature-driven flow, at best it could be called a marginal case where all non-constant terms vanish in the Taylor expansion of $v(\kappa)$.
Nevertheless, it deserves special attention. not only because it is one of the three components of the Bloore flow Eq.(\ref{Bloore}),\cite{Bloore} which is of central interest in this paper,
but also because, due to its simplicity, it will serve in our analytical example. Unlike curvature-driven flows with $v_{\kappa}>0$, the Eikonal flow has no smoothing property,
rather, it can create new singularities and it is driving surfaces \it away \rm from the sphere. The latter can be also seen by regarding the time-reversed (outward moving) Eikonal flow which is converging to the sphere.
Since singularities can arise under the Eikonal flow, here the global evolution of $N(t)$ can only be studied by considering non-smooth shapes as well. We consider only convex shapes, because non-convex, non-smooth shapes can not
be propagated under the Eikonal flow. Our goal is to prove
\begin{lem} \label{lem:Eikonal}
Under the Eikonal flow
N(t) is constant if the curve is  convex and smooth. Under the Eikonal flow with a suitably added truncation rule, specified below,
 $N(t)$ is monotonically decreasing if the curve is convex and piecewise smooth with vertices.
\end{lem}

\begin{figure}[h!]
\begin{center}
\includegraphics[scale=0.5]{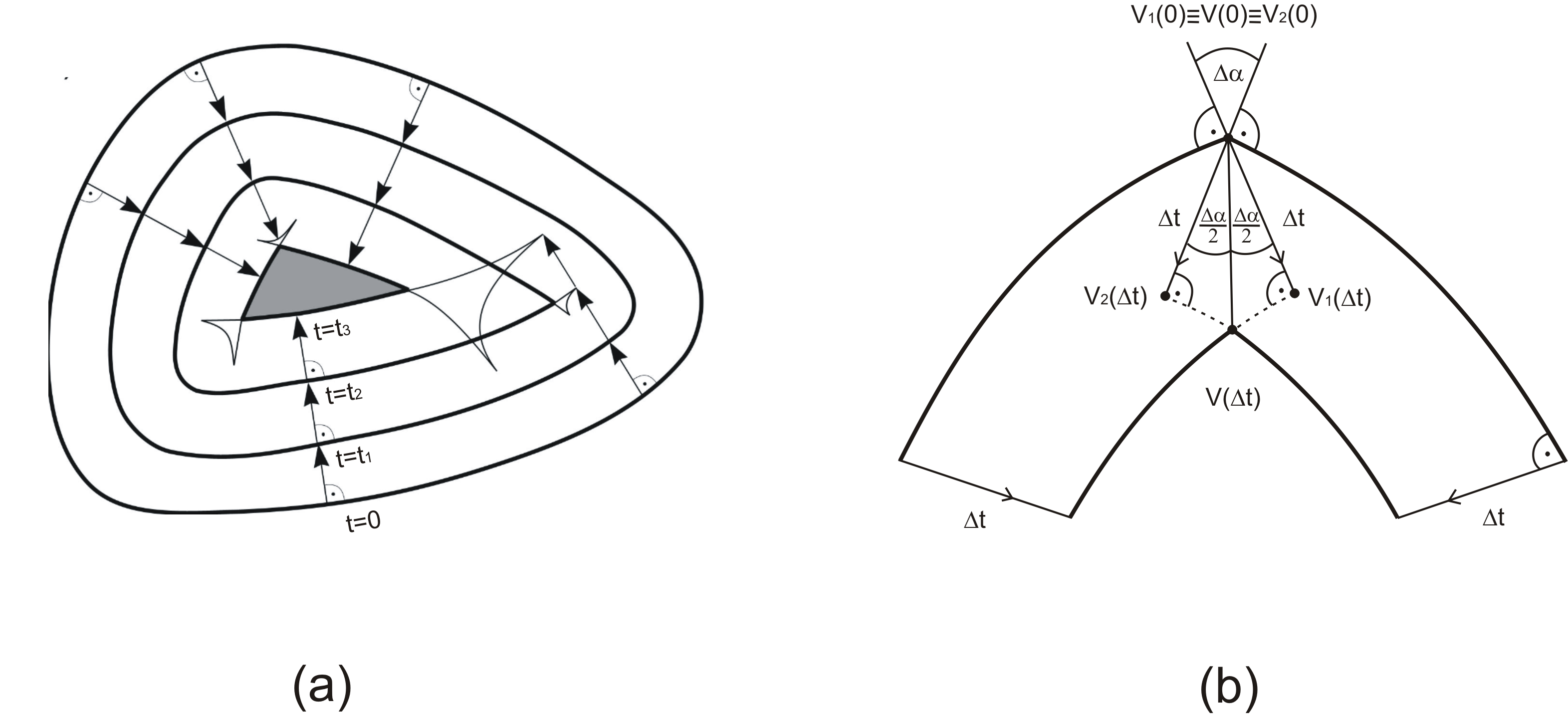}
\end{center}
\caption{\footnotesize{The geometry of  the $v=1$ Eikonal flow Eq.(\ref{Eikonal}). (a) Global geometry of the wavefront starting from a smooth initial condition at $t=0$.
The front remains smooth in the interval $t\in[0,t_1]$ and the first singularity appears in the interval $t \in [t_1,t_2]$.  If the Eikonal flow is serving as an abrasion model
then beyond the singularity the non-convex parts are truncated. (b) Detailed geometry of truncation at a convex vertex. We define the trajectory $V(t)$
of  the vertex as
the bisectrix. Truncated parts of the front marked with dashed line.}}
\label{fig:wavefront}
\end{figure}

\bf Proof. \rm
(a) Smooth, convex curves: As long as the curve is smooth, Theorem 1 predicts that $N(t)$ will remain constant since we have $v_{\kappa}=0$ and thus $\Omega =0$.
One may wonder though whether higher order terms would not cause variation of $N(t)$.  However, under the Eikonal flow, points travel along straight lines
and tangents remain invariant, so, as long as the curve is smooth, we have $N(t)=constant.$ Figure \ref{fig:wavefront} (a) illustrates the global evolution of the wavefront,
which remains smooth for $t\in[0,t_1]$.

(b) Piecewise convex curves: Such curves consist of smooth, convex segments separated by convex vertices $V(t)$ which we will regard at $t=0$ as two coincident
points $V_1(0)$ and $V_2(0)$, belonging to the left and right convex segment, respectively. The relative angle between the normal
at $V_1(0)$ and $V_2(0)$ will be denoted by $\Delta \alpha$. As the front propagates, after time $\Delta t$ the points $V_1(\Delta t)$ and $V_2(\Delta t)$, lying on the corresponding normals,
will be \it separated, \rm simultaneously a self-intersection point, propagating along the bisectrix of $\Delta \alpha$ will
emerge and this point we define as $V(t)$. Since we have no information on how the front should be connected between $V_1(\Delta t)$ and $V_2(\Delta t)$, we \it truncate \rm
the segments $[V_1(\Delta t),V(\Delta t)], [V(\Delta t),V_2(\Delta t)]$. After the truncation we have again a piecewise smooth, convex wavefront.
Based on the argument for the smooth curve we see that the number of critical points on the complete front (without truncation) remains constant. By truncating it, we 
may remove, but never add new critical points, so $N(t)$ will be monotonically decreasing in the non-smooth case under this truncation rule.
Thus we have proven Lemma \ref{lem:Eikonal}.

We also note that a convex vertex is either not a critical point or it is a (non-smooth) maximum. Since minima and maxima can only appear or disappear
in pairs and, as we showed in the first part of the proof, on smooth segments this never happens. The only possible such event is when a smooth minimum reaches the nonsmooth maximum and then they get annihilated due to the truncation, after which the vertex
will become non-critical. As a consequence,
under Assumption 3, substituting the Eikonal equation into Eq.(\ref{alpha5}) yields for the nonsmooth case
\ben \label{Omega_nonsmooth_Eikonal_random}
\bar \Omega= \mbox{sgn}\left(1 +\eta\right).
\een

Finally, we remark that piecewise smooth, convex curves may emerge spontaneously if the Eikonal flow is regarded as an abrasion model \cite{DomokosSiposVarkonyi}.
Unlike the smoothing action of the inward curvature flow, under the Eikonal flow initially smooth curves will develop singularities,
the first will emerge at $t=R_{min}$ where $R_{min}$ is the minimal radius of curvature. In Figure \ref{fig:wavefront} (a) we can see
the spontaneous emergence of singularities in the interval $t_1<t<t_2$. Beyond the singularity, the non-convex parts of the
curve become irrelevant for physical surface evolution models and the internal, convex part will have a piecewise smooth boundary with vertices separated by finite, smooth segments.
We also note that the set of polygons (and that of polyhedra) are invariant under the Eikonal flow.

\subsection{An illustrative example}\label{ss:example}

In case of any fixed reference point (i.e. in the \emph{fixed model}) Lemma \ref{lem:Eikonal} predicts that under the Eikonal flow
$N(t)$ will always decrease on non-smooth convex curves.
Our example is a polygon which is a special case of the latter. 
To quantify our comparison, we introduce the following
\begin{defn}

We call $N(t)$ globally bounded if $\forall t>0$, $N(t) \leq N(0)$.
We call $N(t)$ ultimately bounded if $\exists t_1>0$ such that if $t>t_1$ then $N(t) \leq N(0)$.
\end{defn}
We assume that we can pick initial shapes `uniformly randomly' (which we clarify later) and under this choice
of initial data  we introduce probabilities both in the deterministic and the stochastic model:
\begin{itemize}
	\item in the deterministic model we associate the probabilities $P^{det}_g\geq P^{det}_u$ with globally and ultimately bounded
trajectories, respectively.
\item in the one-parameter ($d$) stochastic model we associate  the probabilities  $P^{st}_g(d)\geq P^{st}_u(d)$ with globally and ultimately bounded
trajectories, respectively.
\end{itemize}
These probabilities admit 
an easy comparison between the models. We can immediately see that due to Theorem 1 and Lemma \ref{lem:Eikonal}, in \it fixed \rm models we have $P^{det}_g= P^{det}_u=1$
and this behavior is imitated in the stochastic model Eq.(\ref{alpha5}) with $d=v_{\kappa}$  and
in case of the non-smooth Eikonal flow Eq.(\ref{Omega_nonsmooth_Eikonal_random}) by setting $d=1$, yielding $P^{st}_g(1)= P^{st}_u(1)=1$.
Our goal is to minimize the error term 
\ben \label{error}
\delta(d) =(P^{det}_u-P^{st}_u(d))^2+(P^{det}_g-P^{st}_g(d))^2,
\een
and  in our example we will show that in case of the fixed model the minimum of $\delta(d)$ is at $d=1$  (corresponding to
the case when no noise is added) and in case
of the centroid model it is at some $d>1$, corresponding to some finite amount of added stochastic noise. This
illustrates that in this example finite stochastic noise added to the fixed model
provides the best approximation to the centroid model, justifying Assumption 3.

\begin{figure}[h!]
\begin{center}
\includegraphics[scale=0.5]{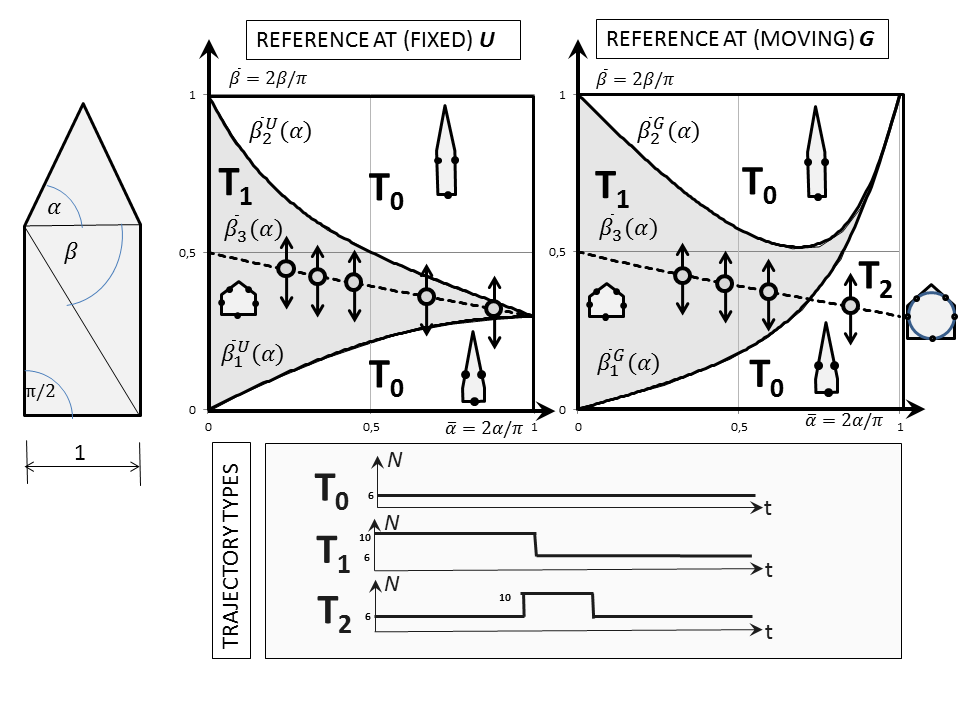}
\end{center}
\caption{\footnotesize{Qualitative evolution of stationary points with respect to the centre of gravity $G$.
Axi-symmetric 5-gon (left),  phase space
with normalized, dimensionless coordinates $\bar \alpha=2\alpha /\pi, \bar \beta= 2\beta /\pi$ on right.
Invariant (repelling) subspace $\bar \beta_3(\alpha)$  marked with dotted line. Parallel flow $v=1$ corresponds to vertical flow in the phase space.
Depending initial data $N(t)$ may remain constant, decrease once, or increase and decrease once subsequently, with $i$ jumps we have $T_i$-type trajectories, marked
in both plots. Time evolution of trajectories shown below. $N(\bar \alpha,\bar \beta)$ has three domains, separated by the curves $\bar \beta_1,\bar \beta_2$.
5-gons in grey middle domain with $\bar \beta_1 < \bar \beta < \bar \beta_2$  have $N=10$
critical points, 5-gons in outer (white) domains have $N=6$ critical points ($S=N/2$ stable points marked on 5gons on insets).}}
\label{fig:angles1}
\end{figure}

Now we introduce our analytical example. Let us regard  the axi-symmetric 5-gon 
in Figure \ref{fig:angles1} given by the angles $\alpha, \beta$ and for simplicity we assume unit width by normalizing the area $A$ as
\ben
A=\frac{1}{4}\tan \alpha + \tan \beta.
\een
The shape is evolving under the Eikonal flow Eq.(\ref{Eikonal}). Since the family of all such shapes
is invariant under Eq.(\ref{Eikonal}), so the PDE is reduced to a system of two first-order ordinary differential equations (ODEs)
\bea
\label{alpha1}
\dot \alpha  & = & f_1(\alpha,\beta)\\
\label{beta1}
\dot \beta & = & f_2(\alpha,\beta).
\eea 
By normalizing $\alpha,\beta$ as $\bar \alpha=2\alpha /\pi, \bar \beta= 2\beta /\pi$, 
the phase space of Eq.(\ref{alpha1})-Eq.(\ref{beta1}) is the unit square, illustrated on the right hand side
of Figure \ref{fig:angles1}.
Since all edges retreat parallel to themselves, in Eq.(\ref{alpha1}) we have 
$f_1=0$ and the flow will be restricted to vertical, $\alpha=constant$ lines and it will be governed by the single ODE (\ref{beta1}).
Self-similar evolution corresponds to fixed points of Eq.(\ref{beta1}) and this can be found directly by observing that if the largest inscribed
circle is tangent to all 5 edges then the 5-gon will evolve in a self-similar fashion. This geometrical condition is
equivalent to $f_2=0$ and it defines an invariant subspace
separating the phase space into two invariant domains (for details see Appendix):
\ben \label{separatrix}
\beta_3= \arctan \left(\frac{1}{2} \left(1+\sqrt{\tan ^2 \alpha+1} -\tan \alpha \right)\right)
\een 
shown as dotted line in Figure \ref{fig:angles1}.
We can also verify that the subspace is repelling, which is indicated by arrows.

So far we have established the global geometry of the phase flow and this is \it independent \rm of the choice of reference point.
Any such choice will determine an integer-valued function $N(\bar \alpha,\bar \beta)$ on the phase space and the flow will then determine
the evolution $N(t)$ along the trajectories of Eq.(\ref{beta1}). We will describe two choices for the reference point:
first we will regard the \it ultimate point  U \rm as a fixed reference and subsequently we will choose the centre of gravity $G$ as a moving reference point.
It is easy to see that in case of the Eikonal flow the centre of the largest inscribed circle
is identical to the ultimate point $U$. We can select $U$ as reference along each $\alpha=constant$ trajectory,
however, it is meaningless to compare fixed reference points along different trajectories. As we can see in Figure 
\ref{fig:angles1}, in both cases $N(\bar \alpha,\bar \beta)$ has three domains, separated by the curves $\bar \beta_1,\bar \beta_2$  (for details see Appendix) and the grey middle domain with $\bar \beta_1 < \bar \beta < \bar \beta_2$  has $N=10$
critical points and the outer domains have $N=6$ critical points (stable points marked on insets). 

The main difference between the fixed and the centroid model is
that in the former case the $N=10$ grey domain contains the entire invariant subspace $\bar \beta_3(\bar \alpha)$ while in the latter case $\bar \beta_3(\bar \alpha)$ and $\bar \beta_1(\bar \alpha)$ intersect transversally
(for details see Appendix). This implies that in the fixed model $O\equiv U$ we can have two types of trajectories:
the ones originating in the grey domain with $N=10$ critical points will decrease once in their evolution  and this we call a type $T_1$
trajectory. Trajectories originating outside this domain will have constant $N(t)$ with no jumps and this we call type $T_0$
trajectories. Using the formulae given in the Appendix, we can compute the corresponding areas as
\ben \label{pd}
q_0^U=0.65, \quad q_1^U=0.35
\een
where superscripts refer to the point of reference. In case of the centroid model with $O\equiv G$
we can distinguish three types of trajectories based on the number of jumps in  $N(t)$. In addition to Type $T_0$ and $T_1$ we can also observe $T_2$-type trajectories with two jumps. Here $N(t)$ will increase once and subsequently decrease also
once by the same amount, so its final value will be equal to the initial value. The corresponding areas can be computed as
\ben \label{qi}
q^G_0=0.49 , \quad q^G_1=0.42, \quad q^G_2=0.09.
\een

\begin{figure}[h!]
\begin{center}
\includegraphics[scale=0.6]{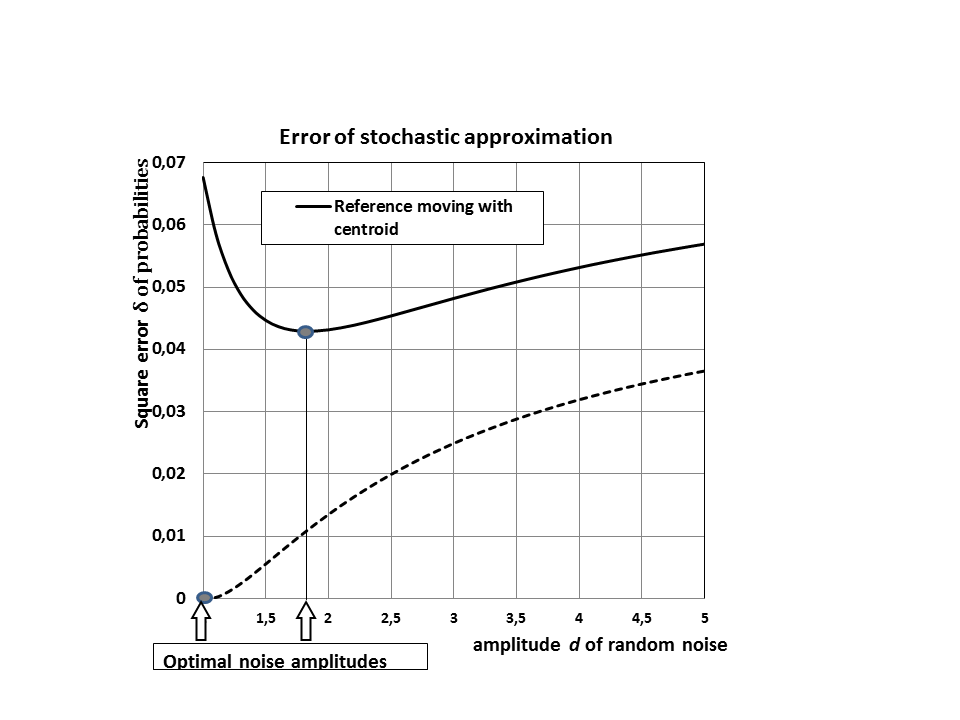}
\end{center}
\caption{\footnotesize{Error of stochastic approximation: square error $\delta$ Eq.(\ref{error}) plotted versus amplitude $d$ of added random noise. Solid line: reference at centre of gravity. Dashed line: reference fixed.
 In case of fixed reference zero noise (deterministic model) is optimal while
in case of reference at centre of gravity random noise initially improves the approximation.}}
\label{fig:errors}
\end{figure}

To obtain a statistical comparison between the models, we now pick initial locations uniformly
randomly on the unit square in the $[\bar \alpha, \bar \beta]$-space and we associate the probabilities $q^U_i, q^G_i$ with the respective trajectory types
and from these values we can obtain the probabilities associated with \it globally bounded \rm
and \it ultimately bounded \rm trajectories, according to Definition 1. As we can immediately see, the fixed model will have unit probability
associated with both categories (since all trajectories are monotonic):
\ben \label{pdetu}
P^{U (det)}_{glob}=P^{U (det)}_{ult}=q_0^U+q_1^U=1.
\een
In case of the centroid model we do see non-monotonic $T_2$ trajectories, so we have
\ben \label{pdetg}
P^{G (det)}_{glob}=q_0^G+q_1^G=0.91, \quad P^{G (det)}_{ult}=q_0^G+q_1^G+q_2^G=1.
\een
We will compare the probabilities associated with these categories in the deterministic model and the stochastic model, in both case we use uniform random choice of initial conditions.
In the stochastic model we pick an initial condition in the described, uniformly random manner and if the corresponding deterministic trajectory is of type $T_i$ then, using Assumption 3
and equation Eq.(\ref{Omega_nonsmooth_Eikonal_random}) we
make $i$ subsequent, independent random draws from the binary distribution
\ben
P(\Delta N=-k)=p,\quad P(\Delta N=+k)=1-p
\een
where $p=(d+1)/2d$. (This is a realization of the Markov chain, introduced in Eq.(\ref{transition})).
Based on these random events we can compute in the stochastic model
\ben \label{stoc}
P^{st}_{glob} =  q_0+q_1p+q_2p^2, \quad
P^{st}_{ult}  =  q_0+2q_2p(1-p)+q_1p+q_2p^2. 
\een
and of course we can do this both for the fixed ($O\equiv U$) and the centroid ($O\equiv G)$ models.
Now we substitute Eq.(\ref{pdetu}),Eq.(\ref{pdetg}) and Eq.(\ref{stoc}) into the error term Eq.(\ref{error}) and minimize the latter by varying the
amplitude $d$ of the random noise. The result is illustrated in Figure \ref{fig:errors} and we can see, that (as expected) in case of the
fixed model  we have zero error at zero noise and the error grows monotonically with $d$. However, in the centroid model the error
has a marked minimum at $d\approx 1.8$ showing that adding random noise is indeed providing a reasonable approximation for the moving reference point.

\section{Conclusions}\label{sec:conclusions}

\subsection{Motivation for the genericity assumption}

Genericity assumptions are, in general, feasible if small perturbations of the investigated functions are admissible.
If we regard the solution set of a differential equation, this is not the case. This indicates that
Assumption 1 on the genericity of bifurcations is a severe restriction from the mathematical point of view, however,
it can be motivated by a heuristic argument.

From the point of view of physical modelling we observe that 
since the PDEs (\ref{PDE1}) and (\ref{PDE_3D}) offer an approximate, averaged continuum description of the discrete event-based
abrasion process, even if the mathematical model predicted degenerate bifurcations, it would be surprising to
actually observe those in the physical process. So, from the physical modelling point of view the genericity
assumption appears to be plausible.

From a more mathematical perspective, by adopting a convenient function norm (e.g. the $C^n$-maximum norm defined
as the supremum over all partial derivatives of order up to $n$, the latter denoting the order
of the first non-vanishing term of the Taylor series) we can see that, on one hand, the differential
operators (\ref{PDE1}) and (\ref{PDE_3D}) are continuous in this norm, on the other hand, functions with only generic
bifurcations (to which we refer as \it generic functions \rm) form a dense, open subset. This implies that the vicinity
of every bifurcation point associated with any
solution  of Eq.(\ref{PDE1}) or Eq.(\ref{PDE_3D}) can be approximated by a convergent sequence of generic functions
and these functions will satisfy the corresponding differential operator with vanishingly small errors.
Since the proof of Theorem 1 depends on a sign (rather than on an exact magnitude), we expect the proof to remain valid if an
arbitrarily small error term is added, implying that sufficiently close generic approximations of solutions have monotonically
decreasing number of critical points. It remains to be shown that among generic approximations there exists always at
least one which is not only in the close vicinity of the approximated solution but also
inherits the number of critical points. If critical points correspond to transversal roots of the first derivative then this
is immediately seen (the first derivative of the generic approximation needs to be a level set of the first derivative
of the solution). In case of multiple roots (arising at tangencies)  we have to select an approximating  generic sequence where the roots are separated
and the separation approaches zero as we converge to the solution. This line of thought suggests that Theorem \ref{thm1}
could be generalized by replacing Assumption 1 by an assumption that smooth curves remain smooth under the flow. 
Verification of the latter assumption may impose further restrictions on $v(\kappa)$.

\subsection{Concluding remarks}

As we pointed out in the Introduction, curvature-driven flows serve both as mathematical tools and as physical models of
surface evolution processes. Our work was motivated by the latter, and in this spirit
 we extended existing results on the evolution of critical points under 
curvature-driven PDEs with the additional assumption of generic bifurcations.

One of the main previous results is due to Grayson \cite{Grayson} on the monotonic decrease of
critical points under the planar curvature flow $v=\kappa$. The other previous result is due to Damon \cite{Damon}
who showed that critical points can be created under the heat equation in three dimensions. The third source was  the result by Kuijper and Florack
 \cite{Kuijper} stating that under a suitable stochastic assumption the number of
critical points is decreasing under the heat equation for surface evolution.
To achieve predictions about curvature-driven flows as models of abrasion,
we combined and partially generalized these ideas. First we showed that under the assumption
of generic codimension-1, saddle-node bifurcations Grayson's result  may be extended to arbitrary 
curvature-driven flows, next we generalized  the result of Kuijper and Florack by weakening the stochastic
assumption and also including arbitrary curvature-driven PDEs.

Our results indicate that in the physical process of mutual abrasion of particles (governed by Eq.(\ref{Bloore}))
the number $N(t)$ of static equilibrium points will decrease 
stochastically, being governed by the asymmetric Markov process Eq.(\ref{transition}). This observation contains no information on the initial and final
values of $N(t)$ in the history of the particle and to predict those is beyond
the scope of the current study. However, we observe that in a natural process, the former is related to
the geometry of fragmentation  while the latter can be observed on pebble beaches among smooth
pebbles, the majority of which tend to be ellipsoidal with just two stable and two unstable equilibrium points.

\section{Acknowledgements}
This work was supported by OTKA grant T104601. The comments and suggestions from  Zsolt G\'asp\'ar, Gary Gibbons, Phil Holmes and two anonymous referees
are gratefully acknowledged. The author is very grateful to Zsolt L\'angi for his invaluable help with computing the derivatives of the principal curvatures with the aid of Maple 16 and for his many other helpful comments.

\bibliographystyle{abbrv} 
\bibliography{References}

\begin{thebibliography}{10}

\bibitem{Alvarez}
L.~Alvarez, F.~Guichard, P.~Lions, and J.~Morel.
\newblock Axioms and fundamental equations of image processing.
\newblock {\em Arch. Rat. Mech. Anal}, 123:199--257, 1993.

\bibitem{Andrews1}
B.~Andrews.
\newblock Contraction of convex hypersurfaces in {E}uclidean space.
\newblock {\em Cal. Var.}, 2:151--171, 1994.

\bibitem{Andrews}
B.~Andrews.
\newblock Gauss curvature flow: the fate of rolling stones.
\newblock {\em Invet. Math}, 138:151--161, 1999.

\bibitem{Andrews2}
B.~Andrews.
\newblock Motion of hypersurfaces by {G}auss curvature.
\newblock {\em Pacific Journal of Mathematics}, 195:1--34, 2000.

\bibitem{Bloore}
F.~J. Bloore.
\newblock The shape of pebbles.
\newblock {\em Math. Geology}, 9:113--122, 1977.

\bibitem{Giblin1}
J.~W. Bruce, P.~J. Giblin, and C.~G. Gibson.
\newblock On caustics of plane curves.
\newblock {\em Amer. Math. Monthly}, 88:651--667, 1981.

\bibitem{Chow1}
B.~Chow.
\newblock On {H}arnack's inequailty and entropy for the {G}aussian curvature
  flow.
\newblock {\em Comm. Pure and Applied Math.}, XLIV:469--483, 1991.

\bibitem{Damon}
J.~Damon.
\newblock Morse theory for solutions to the heat equation and {G}aussian
  blurring.
\newblock {\em J. Diff. Eq.}, 115:368--401, 1995.

\bibitem{jerolmack}
G.~Domokos, D.~Jerolmack, A.~A. Sipos, and A.~T\"or\"ok.
\newblock How river rocks round: explinaining the size-shape paradox.
\newblock {\em PloS One}, DOI:10.1371/journal.pone.0088657, 2014.

\bibitem{Langi}
G.~Domokos and Z.~L\'angi.
\newblock On the equilibria of finely discretized curves and surfaces.
\newblock {\em Monatshefte f\"ur Mathematik}, 168:321--345, 2012.

\bibitem{DomokosSiposVarkonyi}
G.~Domokos, A.~Sipos, G.~Szab\'o, and P.~V\'arkonyi.
\newblock Formation of sharp edges and plane areas of asteroids by polyhedral
  abrasion.
\newblock {\em Astrophysical Journal}, 699:L13--L16, 2009.

\bibitem{Giblin2}
D.~L. Fidal and P.~J. Giblin.
\newblock Generic one-parameter families of caustics in the plane.
\newblock {\em Math. Proc. Camb. Philos. Soc.}, 96:425--432, 1984.

\bibitem{Firey}
W.~Firey.
\newblock The shape of worn stones.
\newblock {\em Mathematika}, 21:1--11, 1974.

\bibitem{Grayson}
M.~A. Grayson.
\newblock The heat equation shrinks embedded plane curves to round points.
\newblock {\em J. Diff. Geom.}, 26:285--314, 1987.

\bibitem{Hamilton}
R.~Hamilton.
\newblock Three-manifolds with positive {R}icci curvature.
\newblock {\em J. Diff. Geom.}, 17:255--306, 1982.

\bibitem{Huisken1}
G.~Huisken.
\newblock Flow by mean curvature of convex surfaces into spheres.
\newblock {\em J. Diff. Geom.}, 20:237--266, 1984.

\bibitem{Huisken}
G.~Huisken.
\newblock Asymptotic behavior for singularities of the mean curvature flow.
\newblock {\em J. Diff. Geom.}, 31:285--299, 1990.

\bibitem{KPZ}
M.~Kardar, G.~Parisi, and Y.-C. Zhang.
\newblock Dynamic scaling of growing interfaces,.
\newblock {\em Phys. Rev. Letters}, 56:889--892, 1986.

\bibitem{kimia1}
B.~Kimia and K.~Siddiqi.
\newblock Geometric heat equation and nonlinear diffusion of shapes and images.
\newblock {\em Computer Vision and Image Understanding}, 64:305--322, 1996.

\bibitem{Koenderink}
J.~Koenderink.
\newblock The structure of images.
\newblock {\em Biol. Cybern.}, 50:363--370, 1984.

\bibitem{Kuijper}
A.~Kuijper and L.~Florack.
\newblock The relevance of non-generic events in scale space models.
\newblock {\em Int.J. of Computer Vision}, 57:67--84, 2004.

\bibitem{Lanczos}
C.~Lanczos.
\newblock {\em Linear Differential Operators}.
\newblock SIAM, Philadelphia, 1996.

\bibitem{Leicht}
K.~Leichtweiss.
\newblock Remarks on affine evolutions.
\newblock {\em Abh. Math. Sem. Univ. Hamburg}, 66:355--376, 1996.

\bibitem{Loog}
M.~Loog, J.~Duistermaat, and L.~Florack.
\newblock On the behavior of spatial critical points under {G}aussian blurring:
  a {F}olklore {T}heorem and scale-space constraints.
\newblock {\em Scale-Space and Morphology in Computer Vision, Lecture Notes in
  Computer Science}, 2106:183--192, 2001.

\bibitem{Mumford}
C.~Lu, Y.~Cao, and D.~Mumford.
\newblock Surface evolution under curvature flows.
\newblock {\em J. Visual Communication and Imgae Representation}, 13:65--81,
  2002.

\bibitem{Maritan}
A.~Maritan, F.~Toigo, J.~Koplik, and J.~Banavar.
\newblock Dynamics of growing interfaces.
\newblock {\em Phys. Rev. Letters}, 69:3193--3195, 1992.

\bibitem{Marsilli}
M.~Marsilli, A.~Maritan, F.~Toigo, and J.~Banavar.
\newblock Stochastic growth equations and reparameterization invariance.
\newblock {\em Rev. Mod. Phys}, 68:963--983, 1996.

\bibitem{Mokhtarian}
F.~Mokhtarian, S.~Abbasi, and J.~Kittler.
\newblock Efficient and robust retrieval by shape content through curvature
  scale space.
\newblock In {\em International Workshop on Image Databases and Multimedia
  Search}, pages 35--42, 1996.

\bibitem{Mokhtarian1}
F.~Mokhtarian and R.~Suomela.
\newblock Robust image corner detection through curvature scale space.
\newblock {\em Pattern Analysis and Machine Intelligence, IEEE Transactions
  on}, 20(12):1376--1381, Dec 1998.

\bibitem{Peitgen}
H.~Peitgen, D.~Saupe, and K.~Schmitt.
\newblock Nonlinear elliptic boundary value problems versus finite difference
  approximations: numerically irrelevant solutions.
\newblock {\em J. Reine u. Angew. Math (Crelle)}, 322:74--117, 1981.

\bibitem{Perelman}
G.~Perelman.
\newblock Ricci flow with surgery on three-manifolds.
\newblock {\em http://arXiv.org/math.DG/0303109v1}, 2003.

\bibitem{Poston}
T.~Poston and I.~Stewart.
\newblock {\em Catastrophe theory and its applications}.
\newblock Pitman, London, 1978.

\bibitem{Rayleigh1}
L.~Rayleigh.
\newblock The ultimate shape of pebbles, natural and artificial.
\newblock {\em Proc. R. Soc. London A}, 181:107--118, 1942.

\bibitem{Rayleigh2}
L.~Rayleigh.
\newblock Pebbles, natural and artificial. {T}heir shape under various
  conditions of abrasion.
\newblock {\em Proc. R. Soc. London A}, 182:321--334, 1944.

\bibitem{Rayleigh3}
L.~Rayleigh.
\newblock Pebbles of regular shape and their production in experiment.
\newblock {\em Nature}, 154:161--171, 1944.

\bibitem{Turyn}
L.~Turyn.
\newblock Spatial critical points of solutions of a one-dimensional nonlinear
  parabolic problem.
\newblock {\em Proc. AMS}, 106:1003--1009, 1989.

\bibitem{Whiteside}
D.~Whiteside, editor.
\newblock {\em The Mathematical Papers of {I}saac {N}ewton, {V}ol 3}.
\newblock Cambridge University Press, 1969.

\bibitem{Oxford}
E.~Zeidler, editor.
\newblock {\em Oxford's User Guide to Mathematics.p 772}.
\newblock Oxford University Press, Oxford, New York, 2004.

\end{thebibliography}

\appendix{
\section{Appendix}

First compute the three curves $\beta^G_1(\alpha), \beta^G_2(\alpha), \beta_3(\alpha)$  in Figure 1.
We use the notations of Figure 5.
The curves $\beta^G_1(\alpha), \beta^G_2(\alpha)$ separate the domains with $N=10$ and $N=6$ critical points.
We can compute these lines based on the conditions that the centre of gravity $G$ should coincide with $P_1,P_2$, respectively. 
In the case of $G\equiv P_1$ we have to write moment balance to the horizontal line passing through $P_1$:
\ben
\frac{h^2}{12}=\frac{a^2}{2}
\een
and by substituting $h=0.5\tan(\alpha)$, $a=\tan(\beta)$ we get
\ben \label{b1}
\beta^G_1(\alpha)=\arctan\left(\frac{\tan(\alpha)}{2\sqrt{3}}\right).
\een

\begin{figure}[!ht]
\begin{center}
\includegraphics[scale=0.5]{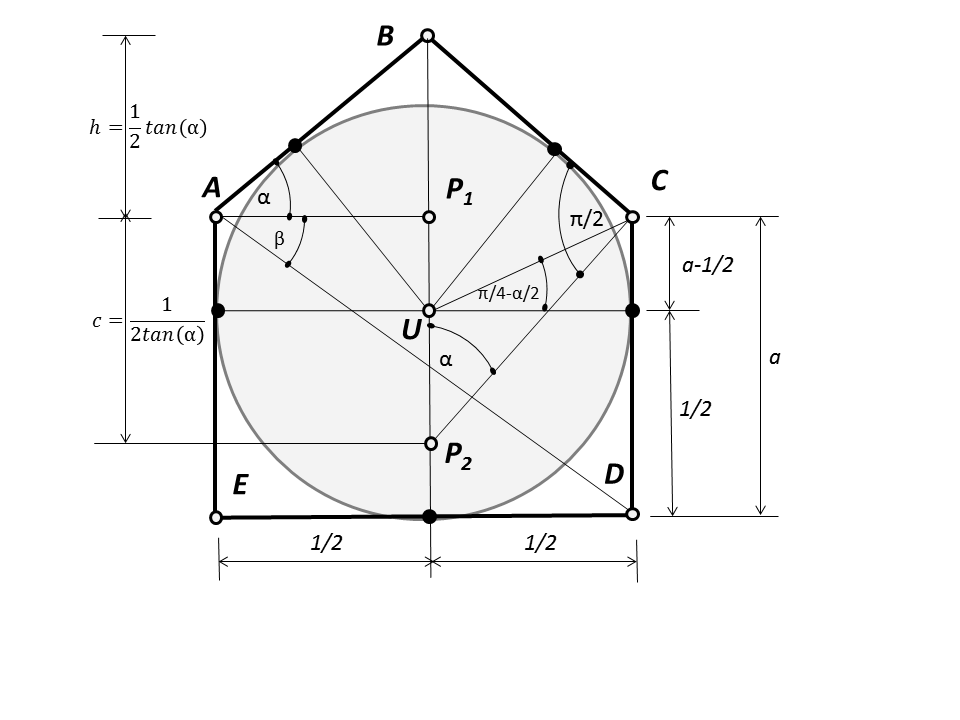}
\end{center}
\caption{\footnotesize{Geometry of the symmetric 5-gon}}
\label{fig:geom}
\end{figure}
In the case of $G\equiv P_2$ we have
\ben
\frac{h}{2}\left(\frac{h}{3}+c\right)=a\left(\frac{a}{2}-c\right)
\een
and by substituting $h=0.5\tan(\alpha)$, $a=\tan(\beta), c=1/(2\tan(\alpha))$ we get
\ben
\frac{\tan(\alpha)}{4}\left(\frac{\tan(\alpha)}{6}+\frac{1}{2\tan(\alpha)}\right)=\tan (\beta)\left(\frac{\tan(\beta)}{2}-\frac{1}{2\tan(\alpha)}\right)
\een
which yields
\ben \label{b2}
\beta^G_2(\alpha)=\arctan \left( \frac{1}{2\tan(\alpha)}\left(1+\sqrt{1 +\frac{\tan ^4(\alpha)}{3}+\tan ^2(\alpha)} \right)\right)
\een
For the invariant subspace we consider the geometry of the 5-gon where all edges are tangent to
the largest inscribed circle and we can write
\ben
a-\frac{1}{2}=\frac{1}{2}\tan\left(\frac{\pi}{4}-\frac{\alpha}{2}\right)
\een
yielding
\ben \label{b3}
\beta_3(\alpha)=\arctan\left(\frac{1}{2}\tan\left(\frac{\pi}{4}-\frac{\alpha}{2}\right)+\frac{1}{2}\right).
\een
In the case of fixed reference point we choose the centre of the largest inscribed circle which is the ultimate point under Eikonal abrasion.
To find the critical curves $\beta^{U}_i(\alpha),(i=1,2)$ we write the conditions for $U\equiv P_i$, yielding
\bea
\beta^{U}_1(\alpha) & = &\arctan\left(\frac{\sin(\alpha)}{2}\right) \\
\beta^{U}_2(\alpha) & = &\arctan\left(\frac{1}{2}\frac{\tan(\alpha)+1}{\tan (\alpha)}\right).
\eea }

\end{document}